\documentclass{article}

\usepackage{amssymb}
\usepackage{amsmath}
\usepackage{amsthm}

\newtheorem{theorem}{Theorem} \newtheorem{lemma}{Lemma}[section]
\newtheorem{propo}{Proposition}[section]
\newtheorem{question}{Question}[section]
\newtheorem{corol}{Corollary}[section] \newtheorem{defin}{Definition}[section]
\newtheorem{remark}{Remark}[section] \newtheorem{claim}{Claim}[section]
\newtheorem{conjecture}{Conjecture}[section] \def\fix {{\#_{\rm fix}}}

 \newcommand{\be}{\beta} \newcommand{\ep}{\varepsilon}
 \newcommand{\de}{\delta} \newcommand{\al}{\alpha}
 \newcommand{\e}{\ep} \newcommand{\te}{\theta}
\newcommand{\HH}{\mathbb{H}} \newcommand{\N}{\mathbb{N}}
 \newcommand{\Z}{\mathbb{Z}}
 \newcommand{\C}{\mathbb{C}}
\newcommand{\F}{\mathbb{F}} \def\finf{\F_\infty}

\def\proof{\smallskip\noindent{\it Proof.} } \def\cur{\curvearrowright}
\newcommand{\OI}{\{0,1\}} \newcommand{\GG}{\mathcal{G}}
\newcommand{\UU}{\mathcal{G}} \newcommand{\G}{\Gamma}
\newcommand{\FFT}{\mathcal{F}_\te} \newcommand{\RR} {\mathcal{R}}
\newcommand{\NR}{\mathcal{N}(\RR)} \newcommand{\ie}{i.\ e.\ }
\newcommand{\hmu}{\hat{\mu}} \newcommand{\trn}{tr_{\mathcal{N}(\RR)}}
\newcommand{\limu}{\lim_{\UU}} \newcommand{\CDR}{\C_d\RR}
\newcommand{\tmd}{tr_{Mat_{d\times d}\mathcal{N}(\RR)}}
\newcommand{\md}{Mat_{d\times d}\mathcal{N}(\RR)}
\newcommand{\mc}{Mat_{d\times d}(\C)}
\title{Sofic equivalence relations\footnote{AMS
Subject Classification: 37A20
\, Research sponsored by OTKA Grant No. 67867, No. 69062}}

\author{G\'abor Elek, G\'abor Lippner}

\begin{document}
\maketitle

\begin{abstract}
  We introduce the notion of sofic measurable equivalence relations. Using
  them we prove that Connes' Embedding Conjecture as well as the Measurable
  Determinant Conjecture of L\"uck, Sauer and Wegner hold for treeable
  equivalence relations.

\end{abstract}

\section{Introduction}
\subsection{Sofic groups and sofic relations}
First let us recall the definition of sofic groups.  The group $\Gamma$ is
sofic if for any real number $0<\epsilon<1$ and any finite subset $F\subseteq
\Gamma$ there exists a natural number $n$ and a function $\psi_n:\Gamma\to
S_n$ from $\Gamma$ into the group of permutations on $n$ elements with the
following properties:
\begin{enumerate}
\item[(a)]   $\fix\Big(\phi(e)\phi(f)\,\phi(ef)^{-1}\Big)\ge(1-\varepsilon)n$
             \kern20pt for any two elements $e,f\in F$.
\item[(b)]   $\phi(1)=1$.
\item[(c)]   $\fix\phi(e)\le\varepsilon n$\kern20pt
             for any $1\ne e\in F$,
\end{enumerate}
where $\fix\pi$ denotes the number of fixed points of the permutation $\pi\in
S_n$.  The notion of soficity was introduced by Gromov \cite{Gro} and Weiss
\cite{Wei} as a common generalization of amenability and residual finiteness.
Direct products, subgroups, free products, inverse and direct limits of sofic
groups are sofic as well. If $N\lhd \Gamma$, $N$ is sofic and $\Gamma/N$ is
amenable, then $\Gamma$ is also sofic.  Residually amenable groups are sofic,
however there exist finitely generated non-residually amenable sofic groups as
well \cite{ESZ}.  It is conjectured that there are non-sofic groups, but no
example is known yet(see also the survey of Pestov \cite{Pes}).

\noindent
In our paper we introduce the notion of a sofic measurable equivalence
relation ($SER$). First let us briefly recall some basic definitions from
\cite{Kech}.  A countable Borel-equivalence relation is a Borel-subspace
$E\subset X\times X$, where $E$ is an equivalence relation and all equivalence
classes are countable. The space $X$ is a standard Borel-space.  Let $\Gamma$
be a countable group and $\Gamma\cur X$ be a Borel-action of $\Gamma$ then it
defines a countable Borel-equivalence relation of $X$ and in fact by the
theorem of Feldman and Moore any countable Borel-equivalence relation can be
obtained by such an action.  A probability measure $\mu$ is $E$-invariant if
it is invariant under a (and actually under all) Borel-action of a countable
group defining the relation $E$.

\noindent
From now on, let $X = \OI^{\N}$ denote the standard Borel space which we equip
with the standard product probability measure $\mu$. For any word $w \in
\OI^k$ $A_w \subset X$ is the closed-open set of those points in $X$ which
start with $w$.  Let $\finf= <\gamma_1,\gamma_2,\dots>$ denote the free group
on countable generators. For any integer $r>0$ let us denoted by $W_r \subset
\finf$ the subset of reduced words of length at most $r$ containing only
letters $\gamma_1, \gamma_1^{-1}, \gamma_2, \gamma_2^{-1} \dots
\gamma_r,\gamma_r^{-1}$. Clearly, $W_0\subset W_1 \subset W_2\dots$ and
$\cup^\infty_{r=0} W_r=\finf$. Suppose $\te : \finf\cur X$ is a (not
necessarily free) Borel group action. Then $\te$ gives rise to a directed
\textit{graphing} (a directed Borel-graph) $\GG \subset X \times X$ in a
natural way: $(x,y) \in \GG$ if and only if there is an index $i$ such that
$\te(\gamma_i,x) = y$. The group action also gives an edge-coloring of this
graphing with countable colors such a way that any vertex there is exactly one
out-edge and one in-edge of every color. The colors are $\gamma_1,
\gamma_1^{-1}, \gamma_2, \gamma_2^{-1} \dots$. Since an edge $xy$ might be
realized by more than one generator, it will be more convenient to think of
$\GG$ as a multi-graphing (\ie one where multiple edges and loop edges are
allowed) and then the action gives us indeed a unique edge-coloring. Also, if
$xy$ is colored by $\gamma_i$ then $yx$ is colored by $\gamma_i^{-1}$.

\begin{defin}
  By an $r$-\textit{neighborhood} we mean an $r$-edge-colored oriented\\
  multi-graph. That is the out-edges need to have different colors from the
  set $\gamma_1, \gamma_1^{-1}, \gamma_2, \gamma_2^{-1} \dots
  \gamma_r,\gamma_r^{-1}$ and if $xy$ is colored by $\gamma_i$ then $yx$ is
  colored by $\gamma_i^{-1}$. Also, we have a chosen vertex which is called
  the root such that any vertex is connected to the root via a path of length
  at most $r$. It is obvious that up to colored, rooted isomorphisms there are
  only finitely many different $r$-neighborhoods. The set of these will be
  denoted by $U^r$.

  Given the group action $\te$ and a point $x\in X$ we define its
  $r$-neighborhood $B_r(x)$ to be the subgraph of $\GG$ spanned by
  $\te(W_r,x)$. Its root is $x$ and it inherits the edge-coloring from $\GG$.
\end{defin}

\begin{defin} By a \textit{$r$-labeled $r$-neighborhood} we mean a
  $r$-neighborhood whose vertices are labeled with words taken from $\OI^r$.
  Again the isomorphism types of such objects form a finite set which we
  denote by $U^{r,r}$.

  Given the group action $\te$ and a point $x \in X$ we define its $r$-labeled
  $r$-neighborhood $B_r^r(x)$ to be the $r$-neighborhood of $x$ with labeling
  defined in the following way: any vertex $y \in B_r(x)$ corresponds to a
  point $y' \in X$. The label of $y$ shall be the unique word $w \in \OI^r$
  for which $y' \in A_w \subset X$.
\end{defin}

For a fixed action $\te$ and a fix $\al \in U^{r,r}$ it is easy to see that
the set $T(\te,\al) = \{x \in X: B_r^r(x) \equiv \al\}$ forms a Borel subset
of $X$. Hence we can take its measure $p_{\al}(\te) = \mu(T(\te,\al))$ which
is clearly a number between 0 and 1.

\noindent
We can repeat everything for any action $\te$ of $\finf$ on a finite set $Y$
whose elements are labeled with elements from $\OI^{\N}$. Then $p_{\al}(\te)$
is defined as $\frac {|T(\te,\al)|}{|Y|}$. We call such vertex labelled sets
$X$-sets.

\begin{defin}
  We say that the Borel action $\te$ is \textit{sofic} if there is a sequence
  of actions $\te_n$ of $\finf$ on finite $X$-sets $Y_n$ such that for any
  $r\geq 1$ and $\alpha\in U^{r,r}$ $\lim_{n\to \infty}
  p_{\al}(\te_n)=p_{\al}(\te)$.
\end{defin}
Note this definition is strongly related to the various notions of graph
convergence (see e.g. \cite{ElekRSA}).

\begin{remark}~\label{finiterem} An action $\te$ is sofic if and only if
  $\te^r = \te|_{\gamma_1,\dots,\gamma_r}$, its restricition to the first $r$
  generators is sofic. The if part follows from choosing a suitable diagonal
  sequence from the sequences $\te^r_n$ that prove the soficity of each
  $\te^r$. For the only-if part one takes the sofic sequence $\te_n$ and
  restricts it to the first $r$ generators, thereby obtaining a sequence
  $\te_n^r$ that is obviously sofic for $\te^r$.
\end{remark} 

\noindent
We call a countable measured Borel-equivalence relation sofic equivalence
relation ($SER$) if it is defined by a sofic action of $\finf$. Obviously,
since any countable group is a quotient of $\finf$, Borel-equivalence
relations can always be defined by $\finf$-actions. In Section
\ref{theoremssec} we shall see that if $E$ is given by actions $\te$ resp.
$\te'$ and $\te$ is sofic, then $\te'$ is sofic 
as well (Theorem \ref{orbeqthm}). That is soficity is not only a property of 
groups actions, but the property of measurable equivalence relations.

\subsection{Results}
We shall prove that Connes' Embedding Conjecture holds for the von Neumann
algebra of a sofic equivalence relation (Theorem \ref{connes}).
Also, any sofic relation satisfies
the Measure-Theoretic Determinant Conjecture of L\"uck, Sauer and Wegner
(Theorem \ref{luck}).  We
also show that \textit{treeable} equivalence relations are always sofic
(Theorem
\ref{tree}). Hence
we prove that the two conjectures above hold for free actions of free groups.

\section{Orbit equivalence}~\label{theoremssec} 

\begin{theorem}~\label{orbeqthm} If $\te_1$ is a sofic action and $\te_2$ is
  measured orbit equivalent to $\te_1$ then $\te_2$ is also sofic.
\end{theorem}

\proof By Remark~\ref{finiterem} it is enough to prove the statement in the
special case when $\te_2$ is obtained from $\te_1$ by adding a generator of
the free group whose action does not change the orbit structure of the
relation. Indeed, from this statement the general case follows easily: to see
that the restriction $\te_2^r$ is sofic add the first $r$ generators of
$\te_2$ to $\te_1$, then restrict to the set of $r$ new generators.

Let $\gamma_1,\dots,\gamma_d,\dots$ generate $\te_1$ and let $\gamma$ denote
the new generator in $\te_2$. Since $\gamma$ does not change the orbit
structure we can find for any point $x \in X$ words $w_x,w'_x \in
<\gamma_1,\dots,\gamma_d,\dots>$ such that $\te_2(\gamma,x) = \te_1(w_x,x)$
and $\te_2(\gamma^{-1},x) = \te_1(w'_x,x)$. In fact we can do this in a Borel
way by taking the shortest and lexicographically smallest $w_x,w'_x$ of all
possible choices.

Let us fix an $\ep > 0$. For this $\ep$ we can find an integer $L$ such that
$\mu(X_0) < \ep/2$ where $X_0 = \{x \in X: |w_x| > L \mbox{ or } |w'_x| > L
\mbox{ or either $w_x$ or $w'_x$ contains a generator $\gamma_i$ where $i >
  L$}\}$. Let us look at $X \setminus X_0$. It is partitioned into a finite
number of Borel subsets $H_i: 1\leq i \leq K$ on which $w_x$ and $w'_x$ are
constant functions of $x$. We shall define a sequence of Borel subsets $X_i
\subset X$ in a recursive way. We start with $X_0$. Then we take $H_1$ and
approximate it by a finite union of standard closed-open subsets of $X$
denoted by $H'_1$ so that $\mu(H_1 \triangle H'_1) \leq \ep/4$. (The
$\triangle$ denotes symmetric difference.) Now let $X_1 = X_0 \cup (H_1
\triangle H'_1)$ and $H''_1 = H_1 \cap H'_1$. Next we take $H_2 \setminus
X_1$, and approximate it by a $H'_2$ which is again a finite union of standard
closed-open subsets of $X$ so that $\mu((H_2 \setminus X_1) \triangle H'_2)
\leq \ep/8$, and set $X_2 = X_1 \cup (H_2 \setminus X_1) \triangle H'_2)$ and
$H''_2 = (H_2 \setminus X_1) \cap H'_2$.  We continue this process for all
$H_i$'s. At each step $H_i \setminus X_{i-1}$ is completely disjoint from each
$H'_j: j < i$ so we can always choose $H'_i$ to be disjoint from all $H'_j: j
< i$. So at the end we have a partition $X = X_K \cup H''_1 \cup \dots \cup
H''_K$ such that $\mu(X_K) \leq \ep, H''_i \subset H_i \cap H'_i$. During the
whole process we considered some large, but finite number of standard
closed-open sets. Each such set is defined by fixing the first few digits of
$x$. Let $M \geq L$ denote an integer such that none of the used closed-open
sets require fixing more than $M$ digits of $x$. Now if $x \in X \setminus
X_K$ then the first $M$ digits of $x$ determine which $H'_i$ it is in, and
hence which $H''_i$ and which $H_i$ it is in. This in turn determines $w_x$
and $w'_x$.

So in fact we have a Borel splitting $X = X_K \cup X'$ such that $\mu(X_K) <
\ep$ and for any point $x \in X'$ the words $w_x,w'_x$ are determined by the
first $M$ digits of $x$.

We have the sofic sequence $G_n$ for $\te_1$. From it we shall construct a
sequence $G_n^{\ep}$. As a first attempt for each vertex $g \in G_n$ we read
the first $M$ digits of its label. Then find the corresponding words $w_x,
w'_x$ we defined above, and trace these words in $G_n$ starting from $g$. If
they end at $h$ and $h'$ respectively then we connect $g$ to $h$ by an
oriented edge labeled $\gamma$ and to $h'$ by an oriented edge labeled
$\gamma^{-1}$. At this point the graph $G_n^{\ep}$ might not be the graph of a
group action: the $\gamma$ edge going from $g$ to $h$ might not be matched by
a $\gamma^{-1}$ edge going from $h$ to $g$. Let us temporarily call such $g$
vertices ``bad''. Let us denote by $\varrho_{\ep}(n)$ the ratio of bad
vertices in $G_n^{\ep}$. By the construction of $G_n^{\ep}$ the badness of a
vertex $g$ is determined by its $(M,M)$ neighborhood in $G_n$. Let us call a
neighborhood $\al \in U^{M,M}(\te_1)$ ``bad'' if its root is a bad vertex.
Hence \[\varrho_{\ep}(n) = \sum_{\al \mbox{\scriptsize{ is bad}}}
p_{\al}(G_n).\] Then if $x \in X$ has neighborhood $\al$ then either $x$ or
$\te_2(\gamma,x)$ has to lie in $X_1$. Hence \[\sum_{\al \mbox{\scriptsize{ is
      bad}}} p_{\al}(\te_1) \leq 2\ep.\] This means that $\limsup_{n \to
  \infty} \varrho_{\ep}(n) \leq 2\ep$. Let us complete the construction of
$G_n^{\ep}$ by keeping the $\gamma$ action for the good vertices, and defining
it arbitrarily for the bad vertices to make it a proper action. This can
always be done: let us denote the set of good vertices by $H$. Then
$\gamma(H)$ is the set of $\gamma$-neighbors of the elements of $H$.
Obviously $|H| = |\gamma(H)|$, and hence $|G_n \setminus H| = |G_n \setminus
\gamma(H)|$. So there is a bijection between these last two sets. This
bijection shall be the action of $\gamma$ and its inverse the action of
$\gamma^{-1}$ on $G_n \setminus H$ and $G_n \setminus \gamma(H)$ respectively.

Let us fix $r$ and a neighborhood $\al \in U^{r,r}(\te_2)$. Let us suppose for
a moment that there are no ``bad'' vertices at all. Then since each $\gamma$
edge is at most an $M$-long path of non-$\gamma$ edges, the $r$-neighborhood
of the $\te_2$ action of any vertex is contained in, and determined by the
$r\cdot M$-neighborhood of the same vertex for the $\te_1$ action. Thus we get
a function $\pi : U^{r\cdot M,r\cdot M}(\te_1) \to U^{r,r}(\te_2)$. Let $B
\pi^{-1}(\al)$. Let $H \subset G_n$ denote those vertices $x \in G_n$ whose
$r$-neighborhood $B_r(x,G_n^{\ep})$ contain a ``bad'' verticex. Then obviously
$x \not \in H$ then $x \in T(G_n^{\ep}, \al)$ if and only if $x \in \cup_{\be
  \in B} T(G_n,\be)$. In other words $T(G_n^{\ep}, \al) \triangle (\cup_{\be
  \in B} T(G_n,\be)) \subset H$. On the other hand if $x\in H$ since
$B_r(x,G_n^{\ep})$ contains the ``bad'' vertex $y$ then also $x \in
B_r(y,G_n^{\ep})$. Hence $H$ is covered by the $r$-neighborhoods of the
``bad'' vertices so $|p_{\al}(G_n^{\ep}) - \sum_{\be \in B} p_{\be}(G_n)| \leq
\varrho_{\ep}(n) \cdot r^r$.

 The same holds for $X$: if
$X_0$ happens to be empty then $p_{\al}(\te_2) = \sum_{\be \in B}
p_{\be}(\te_1)$. However $X_0$ might not be empty, and in this case
$T(\al,\te_2)$ is not necessarily the same as $\cup_{\be \in B} T(\be,\te_1)$.
But if the $r$-neighborhood (by $\te_2$) of a point $x \in X$ is disjoint from
$X_K$, then it cannot belong to the symmetric difference of the two sets
above. Hence
\[ |p_{\al}(\te_2) - p_{\al}(G_n^{\ep})| \leq \sum_{\be \in
  B}|p_{\be}(\te_1) - p_{\be}(G_n)| + (\mu(X_K)+ \varrho_{\ep}(n)) \cdot
  r^r.\] So letting $n \to \infty$ we get that if $\al \in
  U^{r,r}$ then
\[ \limsup_{n\to \infty} |p_{\al}(G_n^{\ep}) - p_{\al}(\te_2)|
\leq 3\ep \cdot r^r.\]

Hence letting $\ep \to 0$ we can choose a suitable
diagonal sequence $G_n'$ from the $G_n^{\ep}$'s to get a sofic
sequence for $\te_2$.

\qed

\begin{corol}~\label{involutionscor}
In the definition of soficity we can take actions of $F_2 * F_2 *
\cdots = F_2^{(*\infty)}$ instead of $F_{\infty}$.
\end{corol}

\proof By Remark~\ref{finiterem} it is sufficient to show this on the level of
finitely generated actions. Let us take an action $\te$ of $F_d$ on $X$ and
consider the underlying simple graphing. It has bounded degree (in fact $2d$
is a bound), hence it can be properly Borel edge-colored by at most
$\tbinom{d^2+1}{2}$ colors (see e.g.~\cite{EL}, section 5.3). Hence the same
equivalence relation can be generated as an action $\te'$ of $F_2^{*d'}$ where
$d' = \tbinom{d^2+1}{2}$. Then according to Theorem~\ref{orbeqthm} $\te$ is
sofic if and only if $\te'$ is sofic.  \qed
\section{The von Neumann algebra of a measurable equivalence relation}

In this section we briefly recall the notion of the von Neumann algebra of an
equivalence relation (\cite{FM}, \cite{LSW}).
Let $\RR\subset X\times X$ be a countable Borel-equivalence relation with an
invariant measure $\mu$. Then one has a natural $\sigma$-finite measure $\hmu$
on the space $\RR$ which is  $\mu$ restricted on $X$ ($X\subset \RR$ is given
by the diagonal embedding).
The groupoid ring of $\RR$; $\C\RR$ is defined as follows.
Let $L^\infty(\RR,\C)$ be the Banach-space of essentially bounded functions
on $\RR$ with respect to $\hmu$. Then
$$\C\RR:=\{K\in L^\infty(\RR,\C)\,\mid\, \mbox{there exists $w_K>0$
such that for almost}$$ $$\mbox{ all  $x\in X$:}
\,K(x,y)\neq 0 \,\,\mbox{or}\,\, K(y,x)\neq 0\,\,\mbox{only for $w_K$
amount of  $y$'s.}\}$$
The $*$-ring structure and a trace is given by:
\begin{itemize}
\item $(K+L)(x,y)=K(x,y) +L(x,y)$
\item $KL(x,y)=\sum_{z\sim x} K(x,z)L(z,y)\,$
\item $K^*(x,y)=\overline{K(y,x)}$
\item $\trn(f)=\int_X K(x,x) d\mu(x)$
\end{itemize}
The von Neumann algebra is constructed by the $GNS$-construction.
The inner product $\langle K,L\rangle=\trn (L^*K)$ defines a pre-Hilbert
structure on $\C\RR$ and by $K\rightarrow KL$ we obtain a representation of
$\C\RR$ on the closure $\HH$ of this pre-Hilbert space. The weak closure of
$\C\RR$ in the operator algebra $B(\HH)$ is the von Neumann algebra $\NR$.
The trace $\trn$ extends to $\NR$ weakly continuously to a finite trace on
$\NR$.

\noindent
In Section \ref{conj} we shall study the matrix ring $\md$ as well. Therefore
in our paper we use the following version of the groupoid ring of $\RR$.
Let $$\CDR:=\{K\in L^\infty(\RR,Mat_{d\times d}(\C))\,\mid\, 
\mbox{there exists $w_K>0$
such that for }$$
 $$\mbox{almost all  $x\in X$:}
 \,K(x,y)\neq 0 \,\mbox{or}\, K(y,x)\neq 0\, \,\mbox{only for $w_K$ amount of
   $y$'s.}\}$$ Then $\CDR$ is isomorphic to $Mat_{d\times d}(\C\RR)$. The
 normalized trace $\tmd(K)$ is defined by
$$\tmd(K):=\int_X\frac{Tr K(x,x)}{d}\,d\mu(x)\,,$$
where $Tr$ is the usual trace on $Mat_{d\times d}(\C)$. Observe that
$\md$ can be obtained via the GNS-construction directly as a weak closure
of $\CDR$.

\section{Approximation theorems} \label{appr}
\subsection{The subalgebra of finite type operators}
Let $\RR$ be a sofic equivalence relation on our standard space $(X,\mu)$
given by a sofic Borel-action $\te:\finf\cur X$. Let $\te_n:\finf\cur Y_n$ be
a sofic approximation as in the Introduction.  We define the subalgebra $\FFT$
(the subalgebra of finite type operators) the following way.  Call an element
$K\in \CDR$ $r$-fine, $K\in \FFT^r$ if for any $\alpha\in U^{r,r}$,
$K(y_1,x_1)=K(y_2,x_2)$ if $x_1,x_2\in T(\te,\alpha)$ and $y_1=wx_1$,
$y_2=wx_2$ for some $w\in W_r$. The following properties are easy to check:
\begin{itemize}
\item $\FFT^1\subset \FFT^2\subset\dots$
\item If $K\in \FFT^r, L\in \FFT^s$ then
$K+L\in \FFT^{\max(r,s)}$, $KL\in \FFT^{r+s}$, $K^*\in\FFT^{2r}$,
$Id\in\FFT^1$.
\end{itemize}
That is $\FFT=\cup^\infty_{r=1}\FFT^r$ is a unital $\star$-subalgebra of
$\CDR$.
\begin{propo}
$\FFT$ is weakly dense in $\CDR$.
\end{propo}
\proof
If $K\in \CDR$ then let $s_K=\sup_{x,y}\|K(x,y)\|$, where $\|\,\|$ is the
usual matrix norm. We say that $\{L_n\}^\infty_{n=1}\subset \CDR$ converge
to $L$ in measure ($L_n\stackrel{\mu}{\to} L$). If :
\begin{itemize}
\item there exist bounds $w$ and $s$ such that for any $n\geq 1$, 
$s_{L_n}\leq s$, $w_{L_n}\leq w$.
\item for any $\e>0$, $\lim_{n\to\infty} \mu(A_\e(n))=0$, where
$$A_\e(n):=\{x\in X\,\mid\, \|L(y,x)-L_n(y,x)|\|>\e\,,\,\mbox{for some $y$}\}$$
\end{itemize}
\begin{lemma}
If $L_n\stackrel{\mu}{\to} L$, then $\{L_n\}^\infty_{n=1}$ weakly converges to
$L$.
\end{lemma}
\proof We need to prove that for any $K\in \CDR$, $\tmd K(L_n-L)\to 0$.
We use the inequality $|\frac{1}{d} Tr(AB)|\leq \|A\|\|B\|\,.$
$$|\tmd K(L_n-L)|=|\int_X\frac{1}{d}\sum_{x\sim z}
Tr(K(x,z)(L_n-L)(z,x)d\mu(x)|\leq $$
$$\leq \int_{A_\e(n)} |\frac{1}{d}\sum_{x\sim z}
Tr(K(x,z)(L_n-L)(z,x)d\mu(x)| +\e w_Ks_K\leq $$
$$\leq \mu(A_\e(n))w_Ks_K(s+s_L) + \e w_Ks_K\,,$$
where $s$ is the bound on the norms of the operators 
$\{L_n-L\}^\infty_{n=1}$. \qed

\vskip0.2in
\noindent
Now for $K\in\CDR$ we construct a sequence in $\FFT$ converging to $K$
in measure.
First let $K'_n\in\CDR$ be defined the following way. Let $K'_n(y,x)=K(y,x)$
if there exists $w\in W_n$ such that $y=wx$, otherwise let
$K'_n(y,x)=0$. Clearly, $K_n\stackrel{\mu}{\to} K$.
Now fix $n\geq 1$. It is easy to see there exist operators $\{K_w\}_{w\in W_n}
\subset \CDR$ such that
\begin{itemize}
\item $K'_n(y,x)=\sum_{w\in W_n} K_w(y,x)$
\item $K_w(y,x)=0$, if $wx\neq y$.
\end{itemize}
Let $f_w(x)=K_w(wx,x)$. Then we have an approximating function $f'_w$ such
that
\begin{itemize}
\item $\mu(x\in X\,\mid \|f_w(x)-f'_w(x)\|>\frac{1}{n})<\frac{1}{n|W_n|}$
\item $f'_w$ is constant on the sets $T(\theta,\alpha)$, if $\alpha\in
  U^{r_w,r_w}$, where $r_w$ is some integer depending on $w$.\end{itemize} Now
let $K_n(y,x)=\sum_{w\in W_n} K'_w(y,x)$, where $K'_w(wx,x)=f'_w(x)$ and
$K'_w(y,x)=0$ if $y\neq wx$. Clearly $K_n\in\FFT$ and $\mu(x\in X\,\mid
\|K_n(y,x)-K(y,x)\|>\frac{1}{n})<\frac{1}{n}\,.$ Therefore
$K_n\stackrel{\mu}{\to} K$. \qed
\subsection{Norm estimates}
Let $A\in\CDR$ and denote by $L_A$ the left-multiplication by $A$ on the
groupoid ring $\CDR$. We give a norm estimate for $L_A$ in terms of
$w_A$ and $s_A$.
\begin{propo} \label{norm}
$\|L_A\|\leq K_dw_A s_A$, where $K_d$ is a constant depending only the
dimension $d$.
\end{propo}
For a matrix $X\in\mc$ $\|M\|_{(d)}$ denote the Frobenius norm, that is
$\frac{Tr(X^*X)}{d}=\|M\|_{(d)}^2$. We have $\|M\|_{(d)}\leq k_d \|M\|$ and 
$\|M\|\leq k_d
\|M\|_{(d)}$
for some constant $k_d$, where $\|M\|$ is the usual matrix norm (the
$l^2$-norm). 
Now let $B\in\CDR$.
Then $\|B\|^2=\tmd(B^*B)=\tmd(BB^*)$ that is
$$\|B\|^2=\int_X\sum_{x\sim y}\frac{Tr B(x,y)B^*(y,x)}{d} d\mu(x)=
\int_X\sum_{x\sim y}\frac{Tr B(x,y)B(x,y)}{d} d\mu(x)=$$
$$=\int_X\sum_{x\sim y} \|B(x,y)\|^2_{(d)}=\int_X t_x \,d\mu(x)\,,$$
where $t_x=\sum_{x\sim z}\|B(x,z)\|^2_{(d)}\,.$
On the other hand, $\|L_AB\|^2=\tmd(B^*A^*AB)=\tmd(A^*ABB^*)$.
Hence, 
$$\|L_AB\|^2=\int_X\sum_{x\sim y}
\frac{Tr A^*A(x,y)B^*B(y,x)}{d} d\mu(x)\leq $$ $$\leq
\int_X\sum_{x\sim y}\| A^*A(x,y)\|_{(d)}\| BB^*(y,x))\|_{(d)} \,.$$
Observe that
$$\| BB^*(y,x))\|_{(d)}=\|\sum_{x\sim z} B(y,z)B(x,z)\|_{(d)}\leq $$ $$\leq
k_d^2\|\sum_{x\sim z} B(y,z)B(x,z)\|\leq k_d^2 \sum_{x\sim z}(\|B(x,z)\|^2+
\|B(y,z)\|^2)\,.$$
Therefore we have the following inequality :
$$\|L_AB\|^2\leq k_d^4s_{A^*A}\int_X\frac{1}{2}\sum_{x\sim y, A^*A(x,y)\neq 0}
(\hat{t}_x+\hat{t}_y)\,d\mu(x)\,,$$
where $\hat{t}_x=\sum_{x\sim z}\|B(x,z)\|^2\,.$
Therefore,
$$\|L_AB\|^2\leq k_d^6 s_{A^*A}w_{A^*A}\|B\|^2\,.$$
Since $w_{A^*A}\leq w_A^2, s_{A^*A}\leq s_A^2$ our proposition follows.
\qed

\vskip0.2in
\noindent
The previous proposition can be applied in the case of finite sets as well.
Let $T$ be a finite set and $K:T\times T\to \mc$ be matrix-valued kernel
function. These kernels form an algebra analogous to $\CDR$.
 Again we can define $s_K:=sup_{x,y}\|K(x,y)\|$ and the width
$w_K$ as the supremal number such that for any $x\in T$, $K_T(x,y)\neq 0$
resp. $K_T(y,x)\neq 0$ for at most $w_K$ $y's$. The normalized trace
$Tr_\star(K)$
is defined as 
$$Tr_\star(K)=\sum_{x\in T}\frac{Tr K(x,x)}{d |T|}\,.$$
Again we have the inner product $\langle K,L \rangle=Tr_\star(L^*K)$ and
$L_A(B)=AB$. The following lemma is the finite version of Proposition
\ref{norm}
\begin{lemma} \label{normlemma}
$\|L_K\|\leq K_d w_K s_K$.
\end{lemma}

Finally, we prove a simple lemma about convergence in measure.
\begin{lemma}
If $L_n\stackrel{\mu}{\to}L$ in $\CDR$ then $\lim_{n\to\infty}
\tmd(L^i_n)=\tmd(L^i_n)$.
\end{lemma}
\proof
The fact that $\lim_{n\to\infty}
\tmd(L_n)=\tmd(L)$ directly follows from the definition.
Since $(L^i_n-L^i)=(L_n^{i-1}-L^{i-1})L_n+L^{i-1}(L_n-L)$
a simple induction implies that $L^i_n\stackrel{\mu}{\to} L^i$\qed
\subsection{Sofic approximation}

 For $K\in\FFT^r$ and $n\geq 1$ let $K_n:Y_n\times Y_n\to \C$ be
defined the following way.
Let $K_n(q,p):=K(y,x)$ if $p\in T(\te_n,\alpha)$, $x\in T(\te,\alpha)$
and $wp=q,wx=y$ for some $w\in W_r$. 
We call $\{K_n\}^\infty_{n=1}$ the sofic approximation of $K$.
\begin{propo}\label{appropo}
Let $K\in\FFT^r, L\in FFT^s$ then
\begin{enumerate}
\item $\|K_n+L_n-(K+L)_n\|_{(n)}\to 0\,.$, where
  $\|A\|_{(n)}= Tr_\star(A^*A)$.
\item $\|K_nL_n-(KL)_n\|_{(n)}\to 0\,.$
\item $\|K_n^*-(K^*)_n)\|_{(n)}\to 0\,.$
\item $Id_n=Id$
\item There exists $C_K>0$ such that $\|K_n\|\leq C_K$, where $\|A\|$ denotes
  the usual norm.
\item $\lim_{n\to\infty}\frac{Tr_*(K^i_n)}{|Y_n|}=\tmd K^i$.
\end{enumerate}
\end{propo}
\proof
We call a sequence $L_n:Y_n\times Y_n\to\mc$ {\it  negligible}
if
\begin{itemize}
\item $\{s_{L_n}\}^\infty_{n=1}$ and $\{w_{L_n}\}^\infty_{n=1}$ are bounded
above.
\item $\lim_{n\to\infty} \frac{|Q_n|}{|Y_n|}=1$, where
$$Q_n=\{x\in Y_n\,\mid\, L_n(x,y)=0, L_n(y,x)=0\,\mbox{for any $y\in
  Y_n$}\}\,.$$
\end{itemize}
It is easy to see that if $\{L_n\}^\infty_{n=1}$ is negligible then
$$lim_{n\to\infty} Tr_{*}(L_n)=0\quad\mbox{and}\, Tr_{*}(L_n^*L_n)=0\,.$$
Observe that $\{K_n+L_n-(K+L)_n\}^\infty_{n=1}$, 
$\{K_nL_n-(KL)_n\}^\infty_{n=1}$ and $\{K_n^*-(K^*)_n\}^\infty_{n=1}$ are
all negligible sequences. Hence $(1.)$,$(2.)$ and $(3.)$ hold. The fourth
statement is trivial and the fifth one immediately follows from
Lemma \ref{normlemma}. 

\noindent
Since $Tr_{*}(K^i_n-(K_n)^i)\to 0$ in order to prove $(6.)$ one only needs to
show that
$$\lim_{n\to\infty} Tr_\star(K_n)=\tmd(K)\,.$$
The right hand side is equal to
$$\sum_{\alpha\in U^{r,r}} \mu(T(\theta,\alpha)) c(K,\alpha)\,,$$
where $c(K,\alpha)=Tr K(x,x)$ if $x\in T(\theta,\alpha)$ and $K\in\FFT^r$.
On the other hand the left hand side of the equation is equal to
$$\sum_{\alpha\in U^{r,r}}\frac{T(\theta_n,\alpha)}{|Y_n|} c(K,\alpha)\,.$$
Thus by the sofic property $(6.)$ follows. \qed
\section{Connes' Embedding Conjecture}
In this section we prove Connes' Embedding Conjecture for the von Neumann
algebras of sofic equivalence relations. First let us very briefly
 recall the conjecture
based on the survey of Pestov \cite{Pes} (see also \cite{Oza}). 
Let $R$ be the hyperfinite
factor. Let $\UU$ be a non-principal ultrafilter on the natural numbers
and $\limu$ be the corresponding ultralimit. Consider the algebra
$B_R\subset \prod^\infty_{n=1} R$, where $\{a_i\}^\infty_{i=1}\in B_R$ iff
$\sum_{i\geq 1} \|a_i\|<\infty\,.$ Let $J\subset B_R$ be the ideal of those 
elements $\{a_i\}^\infty_{i=1}$ such that $\limu Tr_R(a_i^*a_i)=0$, where
$Tr_R$ is the unique finite trace on $R$.
Then $R^\omega:=B_R/J$ is the tracial ultrapower of $R$, a von Neumann algebra
factor with trace 
$$Tr_{\UU}\{[a_i]\}^\infty_{i=1}=\limu Tr_R(a_i)\,.$$
\begin{conjecture} [Connes' Embedding Conjecture]
Every separable factor of type $II_1$ embeds into $R^\omega$.
\end{conjecture}
We confirm the conjecture in the case of von Neumann algebras of sofic
equivalence relations.
\begin{theorem} \label{connes}
Let $\RR$ be a sofic equivalence relation. Then $\NR$ embeds into $R^\omega$.
\end{theorem}
\proof By the result of \cite{PR} it is enough to prove that the weakly dense
$*$-algebra $\FFT$ has a trace preserving $*$-homomorphism into $R^\omega$.
Therefore it is enough to construct (see \cite{Oza})
 unital maps $\psi_n:\FFT\to Mat_{i_n\times
  i_n}(\C)$ for some sequence of integers $\{i_n\}^\infty_{n=1}$ such that
for each $K,L\in\FFT$ the following conditions are satisfied.
\begin{itemize}
\item $\limu \|\psi_n(K)+\psi_n(L)-\psi_n(K+L)\|_{(i_n)}=0.$
\item $\limu \|\psi_n(K)\psi_n(L)-\psi_n(KL)\|_{(i_n)}=0.$
\item $\limu \|\psi_n(K^*)-(\psi_n(K))^*\|_{(i_n)}=0.$
\item $\|\psi_n(K)\|$ is a bounded sequence.
\end{itemize}
Now let $\psi_n(K)=K_n$ as in Section \ref{appr}. Then by Proposition
\ref{appropo} all the conditions above are satisfied. \qed

\section{The Measurable Determinant Conjecture}\label{conj}
The goal of this section is to show that the Measurable Determinant Conjecture
of L\"uck, Sauer and Wegner \cite{LSW} holds for sofic equivalence relations.
Let us recall some basic notions from their paper.
Let $A\in Mat_{d\times d'}(\NR)$. Then $AA^*\in\md$ is a positive,
self-adjoint element. Let $E(\lambda)=\chi_{[0,\lambda]}(AA^*)\in \md$
be the spectral projection corresponding to the interval $[0,\lambda]$ and
$F(\lambda)=\tmd E(\lambda)$ be the associated spectral distribution function.
The Fuglede-Kadison determinant is defined as
$$\det_{\md}(AA^*)=\int_{0^+}^\infty \lambda\,dF(\lambda)\,.$$
The Measurable Determinant Conjecture states that
$$\det_{\md}(AA^*)\geq 1$$ provided that $A\in Mat_{d\times d'}(\Z\RR)$, where
$\Z_d\RR\subset\CDR$ is defined by
$$\Z_d\RR:=\{K\in L^\infty(\RR,Mat_{d\times d}(\Z))\,\mid\, 
\mbox{there exists $w_K>0$
such that for }$$
 $$\mbox{almost all  $x\in X$:}
\,K(x,y)\neq 0 \,\mbox{or}\, K(y,x)\neq 0\,
\,\mbox{only for $w_K$ amount of  $y$'s.}\}$$
\begin{theorem} \label{luck}
If $\RR$ is a sofic equivalence relation, then the measurable determinant
conjecture holds.
\end{theorem}
\proof
First let us suppose that $A$ is an operator of finite type. Then
$AA^*\in\FFT$
and we can consider the sofic approximations $\{A_i\}^\infty_{i=1}$,
$\{A_iA_i^*\}^\infty_{i=1}$. Observe that
\begin{itemize}
\item $\det(A_iA_i^*)\geq 1$. Indeed $A_iA_i^*$ is a a positive matrix with
integer entries (see e.g. the proof of Theorem 3.1 (1) in \cite{LSW}).
\item $\{\|L_{A_iA^*_i}\|\}^\infty_{i=1}$ is uniformly bounded.
\item $\lim_{n\to\infty} Tr_{\star} ((A_iA^*_i)^m)=\tmd ((AA^*)^m)$.
\end{itemize}
Then by Lemma 3.2 of \cite{LSW} $\det_{\md}(AA^*)\geq 1$ holds.

\noindent
Now let $A$ be an arbitrary element and $A_nA_n^*\stackrel{\mu}{\to}AA^*$,
where $\{A_nA_n^*\}^\infty_{n=1}\subset\FFT$. By the previous observation
and Proposition \ref{appropo} the conditions of Lemma 3.2 are satisfied, hence
$\det_{\md}(AA^*)\geq 1$. \qed

\section{Examples of sofic equivalence relations}
\subsection{The Bernoulli shift}

Let $\Gamma$ be a group. We consider the Bernoulli space $\OI^{\Gamma}
= \{f: \Gamma \to \OI\}$ The (right) Bernoulli shift $\te :
\OI^{\Gamma} \times \Gamma \to \OI^{\Gamma}$ is defined by
$\te(f,\gamma_1)(\gamma_2) = f(\gamma_1 \cdot \gamma_2)$.
$\OI^{\Gamma}$ can be identified with $X = \OI^{\N}$ by fixing an
enumeration of $\Gamma: \{\gamma_1,\gamma_2,\dots\}$. Then a $k$-digit
label is just a function $\{\gamma_1,\dots,\gamma_k\} \to \OI$.

\begin{propo} The Bernoulli shift of a sofic group is sofic.
\end{propo}

\proof Let $\Gamma$ be a sofic group generated by $s_1,s_2,\dots \in
\Gamma$.  Any element $\gamma \in \Gamma$ can of course be expressed
as a word in these generators, but this expression is usually not
unique. For later use let us fix for each element $\gamma \in \Gamma$
a word $w_{\gamma}$ that expresses $\gamma$ in terms of the
generators. Let us take a sequence of graphs $G_n$ that prove the
soficity of $\Gamma$.  That is, $G_n$ is a directed graph with each
edge being labeled by some $s_i$ such that each vertex has exactly one
in-edge and one out-edge labeled with each generator. We can also
think of this as a right action of the free group $\finf =
<s_1,s_2,\dots>$ on the vertex set of $G_n$. Furthermore the
neighborhood statistics of $G_n$ converge to that of $\Gamma$'s Cayley
graph on these generators.

We shall label each vertex of $G_n$ with an element of $\OI^{\Gamma}$
so that the labeled neighborhood statistic of $G_n$ will converge to
the labeled neighborhood statistic of $\te$. To do so we first assign
to each vertex of each $G_n$ a random bit. This assignment is simply a
random function $\omega: \cup_{n=1}^{\infty} G_n \to \OI$. Then we
take a vertex $g \in G_n$ and assign to it a function $\omega_g :
\Gamma \to \OI$ by the formula $\omega_g(\gamma) = \omega(g\cdot
w_{\gamma})$. Thus now we have an action $\te_n$ on the
$\OI^{\Gamma}$-labeled space $G_n$. We claim that $p_{\al}(\te_n) \to
p_{\al}(\te)$ for any labeled neighborhood $\al$ for a suitable choice
of $\omega$ (in fact for almost all $\omega$'s).

In order to prove this, we shall first consider $\OI$-labeled
neighborhoods, so let us denote by $V^r$ the set of usual
$r$-neighborhoods where each vertex is labeled with 0 or 1, up to
labeled isomorphism. For an $\al \in V^r$ and a $\OI$-labeled graph
$G$ the notations $T(\al,G)$ and $p_{\al}(G)$ extend naturally.  In
the previous paragraph we described how to obtain a
$\OI^{\Gamma}$-labeling from an $\OI$-labeling for the actions $\te_n$
on $G_n$. It is clear by that construction that the
$U^{r,r}$-neighborhood of a vertex $g$ is determined by the
$V^{r+R}$-neighborhood of the same vertex where $R =
\max_{i=1,2,\dots,r} |w_{\gamma_i}|$.

On the other hand there is a natural $\OI$-labeling on the points of
the Bernoulli-shift: just label each $f: \Gamma \to \OI$ by the value
of $f$ on the identity element. In this way we can talk about the
$V^r$-neighborhoods of points of the Bernoulli-shift, and the
$U^{r,r}$-neighborhoods are again determined by the
$V^{r+R}$-neighborhoods in the exact same fashion. Hence to finish the
proof it is enough to show that $p_{\al}(\te_n) \to p_{\al}(\te)$ for
all $\al \in V^r$ for almost all $\omega$'s.

First let $\al \in V^r$ such that its underlying graph is not
isomorphic to the $r$-neighborhood of the identity of $\Gamma$ in the
Cayley graph. Since the $G_n$ is a sofic sequence for the Cayley
graph, it is immediate that $p_{\al}(\te_n) \to 0$. On the other hand
the Bernoulli-shift is essentially free, hence almost all orbits are
isomorphic to the Cayley graph of $\Gamma$ so $p_{\al}(\te)=0$. 

Now let us consider an $\al \in V^r$ whose graph looks like the Cayley
graph around the identity. We can think that the vertices of $\al$ are
indexed by those elements of $\Gamma$ that have length at most
$r$. Then if $f: \Gamma \to \OI$ is a point in the free part of the
Bernoulli-shift then $f \in T(\te,\al)$ if and only if $f(\gamma) =
\al(\gamma)$ for all elements $|\gamma| < r$. (Here we $\al(\gamma)$
denotes the label written on the vertex of $\al$ corresponding to
$\gamma$.) Hence $p_{\al}(\te) = 1/2^{|\al|}$. All we have to prove
now is

\begin{lemma} For almost all $\omega$'s $p_{\al}(G_n) \to 1/2^{|\al|}$.
\end{lemma}

\proof Let us say that a vertex $g \in G_n$ is normal if its
$r$-neighborhood is isomorphic as a graph to the $r$-neighborhood of
the identity element of the Cayley graph. For any vertex $g\in G_n$
let $X_{g}$ denote a random variable that is 1 if $g \in T(G_n,\al)$
and 0 otherwise. Obviously $P(X_g = 1) = 1/2^{|\al|}$ for any normal
vertex $g$ and 0 otherwise, and
\[p_{\al}(G_n) = \frac{\sum_{g \in G_n} X_g}{|G_n|}.\]
 If all the $X_g$'s were independent, then by the law of large numbers
$p_{\al}(G_n)$ would converge to the limit of its expected value with
probability 1, and this expected value is simply
\[ \lim_{n\to \infty} E(p_{\al}(G_n)) = \lim_{n\to \infty} \sum_{g\in G_n}E(X_g) =
\lim_{n\to \infty} \frac{|\{g \in G_n \mbox{
    normal}\}|}{2^{|\al|}|G_n|} = \frac{1}{2^{|\al|}}.\] The $X_g$'s
    are however not independent, but at least they are independent for
    $g$'s in different graphs, and also $X_{g_1},\dots, X_{g_k}$ are
    jointly independent if $g_1,\dots, g_k \in G_n$ are pairwise far
    from each other, namely $d(g_i,g_j) > r$.

  \begin{lemma} There exists a natural number $l > 0$ (depending on $r$)
    and a partition $\cup_{i=1}^l B_i^n = G_n$ such that if $x \neq y \in
    B^n_i$ then the $r$-neighborhoods of $x$ and $y$ are disjoint.
\end{lemma}
\proof Let $H_n$ be a graph with vertex set $V(G_n)$. Let $(x,y) \in
E(H)_n)$ if and only if $B_{r}(x) \cap B_{r}(y) \neq
\emptyset$. Then $deg(x) \leq r^{r}$ for any $x \in V(H_n)$. Let
$l = r^{r} + 1$ then $H_n$ is vertex-colorable by the colors
$c_1, c_2, \dots , c_l$. Let $B_i^n$ be the vertices coloured by
$c_i$.  \qed

Now for a fix $\ep > 0$ let $B_{i_1}^n,\dots,B_{i_{n_q}}^n$ be those elements
of the partition for which $|B_{i_j}^n| \geq \ep/l$. Then since $\{X_g: g \in
B_{i_j}^n\}$ are jointly independent, by the previous argument we get
\[ \lim \frac{{B_{i_j}^n} \cap T(G_n,\al)}{|B_{i_j}^n|} = \lim E\frac{{B_{i_j}^n} \cap T(G_n,\al)}{|B_{i_j}^n|} = \frac{1}{2^{|\al|}}\] almost
  surely for any choice of $i_j$. An easy calculation now shows that
  setting $B = \cup_{j=1}^{n_q} B_{i_j}^n$ we have
\[ \lim \frac{B \cap T(G_n,\al)}{|B|} = \frac{1}{2^{|\al|}}\] for the same set of $\omega$'s. Since $|G_n
  \setminus B| \leq \ep$, this shows that
\[\frac{1}{2^{|\al|}}-\ep \leq \liminf p_{\al}(G_n) \leq \limsup
p_{\al}(G_n) \leq \frac{1}{2^{|\al|}}+\ep\] almost surely, and finally
letting $\ep \to 0$ we get the desired almost sure convergence.

\qed

Thus we have $p_{\al}(\te_n) \to p_{\al}(\te)$ almost surely for all
$\al$'s. Hence there exists an $\omega$ for which $p_{\al}(\te_n) \to
p_{\al}(\te)$, hence the Bernoulli shift is sofic.

\qed

\vskip 0.2in
\noindent
Note that the fact that for residually amenable groups the Measurable
Determinant Conjecture holds for the Bernoulli shift has already been
proved in 
\cite{LSW}.

\subsection{Treeable relations}
Recall \cite{Kech} that
an equivalence relation $E \subset X \times X$ is called {\it treeable}
if it has an
L-treeing generated by measure-preserving involutions
$S_1,S_2, \dots$. We prove that all treeable equivalence relations
are sofic. The most important examples of such treeable relations
are the free actions of free groups.

\begin{theorem} \label{tree} The action of $\Gamma = <\gamma_1,\gamma_2, \dots|
  \gamma_i^2 = 1 (i=1,2,\dots)>$ defined by $\te(\gamma_i,x) = S_i(x)$
  is sofic.
\end{theorem}

\proof By Remark~\ref{finiterem} it is again sufficient to work with
finitely generated actions. So let us assume $\Gamma$ is generated by
$\gamma_1,\dots,\gamma_d$.  Let us fix a large $r$. For any $\al, \be
\in U^{r,r}$ and any $1 \leq i \leq d$ let us denote
\[ T(\te,\al,i,\be) = \{x \in T(\te,\al): S_i(x) \in T(\te,\be)\}\]
and it measure (as it is obviously a Borel set)
\[ p_{\al i \be}(\te) = \mu(T(\te,\al,i,\be)).\]
There numbers together with the $p_{\al}(\te)$'s satisfy certain equations: 
\begin{eqnarray*}
\sum_{\al \in U^{r,r}} p_{\al}(\te) & = & 1 \\ \sum_{\be \in U^{r,r}}
p_{\al i \be}(\te) & = & p_{\al}(\te) \mbox{ for any $i$} \\ p_{\al i
\be}(\te) & = & p_{\be i \al}(\te) \mbox{ for any $\al, i, \be$.}
\end{eqnarray*}

Let us introduce variables $w_{\al}: \al \in U^{r,r}$ and $w_{\al i
\be} : \al,\be \in U^{r,r}, 1 \leq i \leq d$. Then $w_{\al} =
p_{\al}(\te), w_{\al i \be} = p_{\al i \be}(\te)$ is a solution to the
following set of linear equations:

\begin{eqnarray}
\sum_{\al \in U^{r,r}} w_{\al} & = & 1 \label{1eq}\\ \sum_{\be \in U^{r,r}}
w_{\al i \be} & = & w_{\al} \mbox{ for any $i$} \label{2eq}\\ w_{\al i
\be} & = & w_{\be i \al} \mbox{ for any $\al, i, \be$. \label{3eq}}
\end{eqnarray}

Now we use the rational approximation trick of Bowen \cite{Bowen}.
Let us fix a small $\ep > 0$. If a set of linear equations with
rational coefficients has some solution, then it also has a rational
solution in which each variable is at most $\ep$-far from the
corresponding value of the initial solution. Further we may also
assume that if a variable was 0 in the initial solution then it
remains 0 in the new solution. So our set of equations has such a
rational solution which will shall simply denote by $w_{\al},w_{\al i
\be}$. Since now these numbers are all rational, we may choose a large
integer $N$ for which $W_{\al i \be} = N\cdot w_{\al i \be}$ is always
an even integer.

Now take a set $Y$ with $N$ elements and partition it into subsets
$Y_{\al} : \al \in U^{r,r}$ with $|Y_{\al}| = W_{\al}$. This can be
done because of (\ref{1eq}) above. Then fix an index $i$ and do the
following: if for a type $\al$ the involution $S_i$ is fixing the
root, then define $S_i(y) = y$ for all $y \in Y_{\al}$. Otherwise
partition $Y_{\al}$ into subsets $Y_{\al i \be}$ of size $W_{\al i
\be}$. This can be done because of (\ref{2eq}) above. Finally define
$S_i$ to be a random bijection between $Y_{\al i \be}$ and $Y_{\be i
\al}$, or a random matching in $Y_{\al i \al}$ (this is where we need
that the size of this set is even). This can be done because of
(\ref{3eq}).  Repeat this procedure for each index. Finally for any
$\al \in U^{r,r}$ and any $y \in Y_{\al}$ look at the label of the
root in $\al$. This is a word $w \in \OI^k$. Label $y$ with any
infinite $w' \in \OI^{\infty}$ which starts with $w$.

This way we defined an action $\te'$ of $\Gamma$ on the finite labeled
set $Y$. We claim this will be a good approximation to the action
$\te$. To make this precise let us fix an ordering of all possible
neighborhood types $\al_1,\al_2,\dots$, and for two actions $\te,\te'$
let us introduce their statistical distance $d_s(\te,\te')=
\sum_{i=1}^{\infty} \frac{|p_{\al}(\te) - p_{\al}(\te')|}{2^i}$. It is
easy to see that $\te_n$ is a sofic sequence for $\te$ if and only if
$d_s(\te,\te_n) \to 0$.

\begin{lemma} Let $\nu_q$ denote the ratio of those points in $Y$ through
  which there is a $\te'$ cycle of length at most $q$.  Then for any fixed $q$
  we have $\nu_q \to 0$ in probability when $N \to \infty$.
\end{lemma}

\proof By the construction of $Y$ the probability of the existence of any
particular $xy$ edge is at most $c/N$ for some universal constant $c$
depending only on the $w_{\al i \be}$ numbers. Hence the probability that a
particular cycle of length $l$ exists in $\te'$ is at most $c^l/N^l$, hence
the expected number of length $l$ cycles is at most $c^l/N^l \cdot
\tbinom{N}{l} < c^l/l!$ which is a constant. So for fixed $q$ and large $N$
the expected number of points through which there is cycle of length at most
$q$ is at most some constant $c_q$. Then for any fixed $\ep$ we have $P(\nu_q
> \ep) \leq \frac{c_q}{\ep N}$ so clearly $P(\nu_q > \ep) \to 0$ as $N \to
\infty$. \qed

Then the ratio of those vertices whose $r$-neighborhood is not a tree
is at most $d^r \nu_{2r}$ since any such neighborhood contains a cycle
of length at most $2r$, and hence the root of this neighborhood is at
most $r$ steps from a vertex in the cycle.

For a neighborhood $\al \in U^{r,r}$ let us denote by $\al|_q \in
U^{q,q}$ the subgraph of $\al$ spanned by the vertices that are at
most $q$ steps from the root and keeping only the first $q$ digits of
the labels. The following is easily verified by induction on $q$:

\begin{claim} If $q \leq r$ and the girth of $\te'$ at $y \in Y_{\al}$ is 
greater than $2q$ then $B_q(y) \cong \al|_q$.
\end{claim}

Now we can estimate $d_s(\te,\te')$. Let us fix $r$ and let $j$
denote the index of the first $\al_i$ neighborhood in our listing
either whose radius is larger than $r$ or its labels have more than
$r$ digits. 

Let 
\[ U = \bigcup_{q \leq r} U^{q,q}, \quad U_c = \{ \al \in U:
\al \mbox{ is not a tree}\}, \quad U_t = U \setminus U_c.\] If $\al \in
U_c$ then $p_{\al}(\te) = 0$ since $\te$ is a treeing, and
$p_{\al}(\te') \leq d^r \nu_{2r}$ since at most this many vertices can
have cycles in their $r$-neighborhood.

If $\al \in U_t \cap U^{q,q}$ then 
\begin{multline*} |p_{\al}(\te)-p_{\al}(\te')| = \left| \sum_{\be \in
    U^{r,r}: \be|_q \cong \al} p_{\be}(\te)-p_{\be}(\te')\right| \leq
    \\ \leq \left|\sum_{\be \in U^{r,r}: \be|_q \cong \al}
    p_{\be}(\te)-w_{\be}\right| + d^r\nu_{2r} \leq \ep |U^{r,r}| + d^r\nu_{2r}.
\end{multline*}
The term $d^r\nu_{2r}$ appears again because $p_{\be}(\te')$ is not
necessarily equal to $w_{\be}$: the difference comes from exactly
those vertices in $Y_{\be}$ whose $2r$ neighborhood is not a tree. 
And finally

\begin{multline}
\label{dseq} d_s(\te,\te') = \sum_{i=1}^{\infty} 
\frac{|p_{\al_i}(\te)-p_{\al_i}(\te')|}{2^i} \leq \sum_{i: \al_i \in U_c}
  \frac{d^r\nu_{2r}}{2^i} + \sum_{i: \al_i \in U_t} \frac{\ep
  |U^{r,r}| + d^r\nu_{2r}}{2^i}  + \\ + \sum_{i \geq j} \frac{1}{2^i} \leq
  \ep |U^{r,r}| + 2d^r \nu_{2r} +  1/2^{j-1} 
\end{multline}

So in order to construct a finite action with $d_s(\te,\te') < \de$
first we choose $r$ so large that $1/2^{j-1} < \de/3$ in
(\ref{dseq}). Then we choose an $\ep <
\frac{\de}{3|U^{r,r}|}$. Then we find a rational solution to our
system of equations (\ref{1eq},\ref{2eq},\ref{3eq}). Finally we choose
$N$ so large, that with positive probability $\nu_{2r} \leq
\frac{\de}{6d^r}$. We pick an action $\te'$ satisfying this and hence 
\[d_s(\te,\te') \leq \ep|U^{r,r}| + 2d^r \nu_{2r} + 1/2^{j-1} < 3\cdot 
\de/3 = \de.\]
Hence $\te$ is indeed a sofic action. \qed

\vskip0.2in
\noindent
Note that the previous theorem combined with Theorem \ref{orbeqthm}
shows the all treeable groups are sofic. Recall that a group is treeable
if it has a free treeable action.

\subsection{Profinite actions}
The simplest case of sofic action is argueably the case of profinite actions.
Let $\G$ be a countable residually finite group and
$\Gamma\supset N_1\supset N_2\dots$ be finite index normal subgroups such that
$\cap_{i=1}^\infty N_i=\{1\}$. Then $G=\lim_{\leftarrow} \Gamma/N_i$ is the
profinite closure with respect to the system $\{N_i\}$, a compact group.
Then $\G$ is a dense subgroup of $G$ and so it preserves the Haar-measure
$\nu$. It is easy to see that $\G\cur (G,\nu)$ is a sofic action. 

\section{Conclusion}
We can conclude that the Connes Embedding Conjecture and the Measurable
Determinant Conjecture hold for treeable sofic relations, particularly, for
relations induced by free actions of free groups. We end our paper with a
question related to Question 10.1 of Aldous and Lyons \cite{AL}
on unimodular networks.
\begin{question}
Does there exist a measurable equivalence relation that is not sofic ?
\end{question}


\begin{thebibliography}{9}
\bibitem{AL} {\sc D. Aldous and R. Lyons}, {\em 
Processes on Unimodular Random \\
Networks}, {\sl  Electron. J. Probab.}  12  (2007), no. 54, 1454-1508.
\bibitem{Bowen} {\sc L. Bowen}, Periodicity and circle packings of the 
hyperbolic plane. {\sl  Geom. Dedicata} {\bf  102}  (2003), 213-236. 
\bibitem{ElekRSA} {\sc G. Elek},
A Regularity Lemma for Bounded Degree Graphs and Its Applications: 
Parameter Testing and Infinite Volume Limits (preprint) arXiv:0711.2800
\bibitem{EL} {\sc G. Elek and G. Lippner},
An analogue of the Szemeredi Regularity Lemma for bounded degree graphs,
(preprint) arXiv:0809.2879
\bibitem{ESZ} {\sc G. Elek and E. Szab\'o},
On sofic groups, {\sl J. Group Theory} {\bf 9}  (2006)  no. 2, 161--171.
\bibitem{FM} {\sc J. Feldman and C. C. Moore},
Ergodic equivalence relations, cohomology, and von Neumann algebras. II. 
{\sl Trans. Amer. Math. Soc.} {\bf 234} (1977), no. 2, 325-359.
\bibitem{Gro} {\sc M. Gromov},
Endomorphisms of symbolic algebraic varieties
{\sl J. Eur. Math. Soc.} {\bf 1} (1999) no. 2, 109-197.
\bibitem{Kech} {\sc A. Kechris and B. Miller},
Topics in orbit equivalence. {\sl Lecture Notes in Mathematics} {\bf 1852}
Springer Verlag
\bibitem{LSW} {\sc W. L\"uck, R. Sauer and C. Wegner},
L2-torsion, the measure theoretic determinant conjecture, and 
uniform measure equivalence. (preprint) arXiv: 0903.2925. 
\bibitem{Oza} {\sc N. Ozawa}, 
About the QWEP conjecture. {\sl  Internat. J. Math.} {\bf 15}  (2004), 
 no. 5, 501-530. 
\bibitem{Pes} {\sc V. Pestov}, Hyperlinear and sofic groups: a brief guide. 
{\sl Bull. Symbolic Logic} {\bf  14}  (2008)  no. 4, 449--480. 
\bibitem{PR} {\sc C. Pearcy and J. R. Ringrose},
Trace-preserving isomorphisms in finite operator algebras. {\sl
  Amer. J. Math.} {\bf  90}  (1968) 444-455. 
\bibitem{Wei} {\sc B. Weiss},
Sofic groups and dynamical systems ({\it Ergodic theory and harmonic analysis,
Mumbai, 1999)}
{\sl Sankhya Ser. A.} {\bf 62} (2000) no. 3,  350-359.
\end{thebibliography}
\end{document}